# The "Eating up" Assumption, Transversality Conditions, and Steady States for Discrete Time Infinite Horizon Problems

Dapeng Cai [a] and Takashi Gyoshin Nitta [b]

[a] Institute for Advanced Research, Nagoya University, Furo-cho, Chikusa-ku, Nagoya, 464-8601, Japan; Email: cai@iar.nagoya-u.ac.jp; [b] Department of Mathematics, Faculty of Education, Mie University, Kurimamachiya 1577, Tsu, 514-8507, Japan; Email: nitta@edu.mie-u.ac.jp

**Abstract**

In this paper, we consider how to construct the optimal solutions for a general discrete time infinite horizon optimal control problem. We establish necessary and sufficient conditions for optimality in the sense of a modified optimality criterion. We also consider how transversality conditions are related to steady states. The results are applied to two examples to demonstrate how the new transversality conditions derived in this paper differ from those given in the previous literature.

# 1. Introduction

There has been a large literature that considers the necessity and sufficiency of transversality conditions (TVCs) for infinite-horizon optimization problems with possibly unbounded objectives, which are extremely common in economics, especially in the endogenous growth literature (Ref. [7], [9], [10]). In this paper, we consider how to construct the optimal solutions for a general discrete time infinite horizon optimal control problem. We establish necessary and sufficient conditions for optimality in the sense of a modified optimality criterion. We also consider how transversality conditions are related to steady states. The results are then applied to two examples to demonstrate how the new transversality conditions derived in this paper differ from those given in the previous literature.

# 2. General Definitions and Conventions

For any function $f$, $\text{dom}(f)$ denotes the domain of $f$. Let $f \in \mathbb{R}^{[a,b]}$ be a bounded function. We assume that $f$ is Riemann integrable, and the number $\int_a^b f(t)dt$ is called the Riemann integral of $f:[a,b] \to \mathbb{R}$. Furthermore, if $g \in \mathbb{R}^{[0,\infty]}$ is a bounded function, then we define the improper Riemann integral of $g$ as

$\int_0^\infty g(t)dt \equiv \lim_{b \to \infty} \int_0^b g(t)dt$, provided that the limit on the right-hand side exists.

For any $a \in \mathbb{R}$ and $\{a_t\}_{t=0}^\infty \in \mathbb{R}^\infty$, we write $\limsup a_t = a$ if (i) for any $\varepsilon > 0$, there exists a $T' > 0$ such that $a_t < a + \varepsilon$ for all $t > T'$, (ii) for any $\varepsilon > 0$ and $t \in \mathbb{N}$, there exists an integral $k > t$ such that $a_k > a + \varepsilon$. The expression $\liminf a_t$ is defined dually. Define $b_t$ similarly. We shall use the following properties:

$$\liminf a_t \equiv -\limsup(-a_t),$$

$$\limsup(a_t + b_t) \leq \limsup a_t + \limsup b_t,$$

$$\liminf(a_t + b_t) \geq \liminf a_t + \liminf b_t,$$

$$\liminf(a_t + b_t) \leq \limsup a_t + \liminf b_t \leq \limsup(a_t + b_t).$$

We say that $f : \mathbb{R} \to \mathbb{R}^n$ is differentiable at $x$ iff $\lim_{t \to x} \frac{f(t) - f(x)}{t - x} \in \mathbb{R}^n$, in which case $f'(x)$ (or $\frac{d}{dt} f(x)$) equals precisely to this number. If $f$ is differentiable at each $x$, we simply say that $f$ is differentiable. Similarly, define $f : \mathbb{R}^2 \to \mathbb{R}$, then $f_1(y, z)$ denote the partial derivative of $f$ with respect to $y$; and $f_2(y, z)$ is defined similarly.

Let $\lim_{z \to \infty} f(y, z) = f(y)$, then $f(y, z)$ converged uniformly to $f(y)$ iff for $\forall \varepsilon > 0$, there exists a $Z_\varepsilon$ that is independent of $y$ such that $|f(y, z) - f(y)| < \varepsilon$ whenever $z > Z_\varepsilon$.

We use the following conventions.

"$\forall a \in A : B$" means "for all $a$ in $A$ such that $B$ is true."

"$\forall t \geq 0$" and "$(t)$" are omitted for a binary relation true for all $t \geq 0$.

## 3. The Model and Results

### 3.1 The Setting

We examine the following deterministic discrete time infinite horizon optimal control problem, in which the objective function is allowed to depend on time:

(1)
$$\begin{cases} \text{"}\max \sum_{0}^{\infty} v(x(t), x(t+1), t)\text{"} \text{ subject to} \\ x(0) = x_0, \ \{x(t)\}_{t=0}^{\infty} \in X \subset \mathbb{R}^{\infty}, \end{cases}$$

where $v$ is an extended real-valued function. Notice that the objective functional of (1) is not necessarily finite. It is to be maximized under the optimality criterion to be defined below.

We define an attainable path as follows (Ref. [6]):

**DEFINITION 1.** An infinite sequence $\{x(t)\}_{t=0}^{\infty}$ is said to be attainable if it starts at time 0 from $x(0) \equiv x_0$, with the state and decision being $\{x(t)\}$ such that $\{x(t)\}_{t=0}^{\infty} \subset X$, for all $t \geq 0$.

For simplicity, we only consider reversible investment. Goods are assumed not to be wasted, which is incorporated into the model as a boundary condition. The attainable paths are then converted accordingly: the values of both consumption and capital at a newly chosen "terminal time", $T$, are modified so that both are 0 at $T+1$:

(2) $\quad \forall T, \quad \tilde{x}^T(t) = \begin{cases} x(t), & 0 \leq t \leq T, \\ 0, & t \geq T+1. \end{cases}$

The image of this conversion is as follows: the "doomsday", $T$, suddenly arrives, and the agents are able to "eat up" the entire stocks of goods at short notice before the end of the world and therefore, nothing is wasted. In what follows, we consider a modified overtaking criterion to be defined below, modified to incorporate the boundary condition that goods are not to be wasted.

**DEFINITION 2.** (Modified Overtaking Criterion) $\{x_1(t)\}$ and $\{x_2(t)\}$ are two infinite horizon attainable paths. For all $v \geq 0$, all $t$, $\{x_2(t)\}$ overtakes $\{x_1(t)\}$ if

$$\liminf_{T \to \infty} \left( \sum_{t=0}^{T} \left( v\left( \tilde{x}_2^T(t), \tilde{x}_2^T(t+1) \right), t \right) - \sum_{t=0}^{T} \left( v\left( \tilde{x}_1^T(t), \tilde{x}_1^T(t+1) \right), t \right) \right) > 0.$$

If $\tilde{x}_1^T(t)$ and $\tilde{x}_2^T(t)$ in Definition 2 are replaced with $x_1(t)$ and $x_2(t)$, we derive the overtaking criterion á la Brock (Ref. [1], [2], [5], [11]):

**REMARK 1.** (Overtaking Criterion á la Brock) $\{x_1(t)\}$ and $\{x_2(t)\}$ are two infinite horizon attainable paths. For all $v \geq 0$, all $t$, $\{x_2(t)\}$ overtakes $\{x_1(t)\}$ if

$$\liminf_{T \to \infty} \left( \sum_{t=0}^{T} \left( v\left( x_2(t), x_2(t+1) \right), t \right) - \sum_{t=0}^{T} \left( v\left( x_1(t), x_1(t+1) \right), t \right) \right) > 0.$$

Based on the notion of weak maximality, the optimality criterion is defined as follows.

**DEFINITION 3.** An infinite horizon attainable path $\{x^*(t)\}$ is optimal if no other infinite horizon attainable path $\{x(t)\}$ overtakes it:

$$\liminf_{T \to \infty} \left( \sum_{t=0}^{T} \left( v\left(\tilde{x}^T(t), \tilde{x}^T(t+1)\right), t \right) - \sum_{t=0}^{T} \left( v\left(\tilde{x}^{*T}(t), \tilde{x}^{*T}(t+1)\right), t \right) \right) \leq 0.$$

Note that the above specifications generalize the settings in Ref. [3], [4]. On the other hand, If $\tilde{x}_1^{*T}(t)$ and $\tilde{x}_2^{*T}(t)$ in Definition 3 are replaced with $x_1^*(t)$ and $x_2^*(t)$, we derive the weak maximality á la Brock:

**REMARK 2.** An infinite horizon attainable path $\{x^*(t)\}$ is optimal if no other infinite horizon attainable path $\{x(t)\}$ overtakes it:

$$\liminf_{T \to \infty} \left( \sum_{t=0}^{T} \left( v\left(x(t), x(t+1)\right), t \right) - \sum_{t=0}^{T} \left( v\left(x^*(t), x^*(t+1)\right), t \right) \right) \leq 0.$$

Denote $\{x_T^*(t)\}$ to be the optimal solution for the maximization problem of the finite horizon ($T$) version of (1), subject to (1) and $x(T+1) = 0$; and denote $\{x^\circ(t)\}$ to be its limit as the horizon $T$ grows to infinity, i.e., $x^\circ(t) \equiv \lim_{T \to \infty} x_T(t)$.

The following theorems generalize Theorem 1 of Ref. [3], which only considers stationary utility and production functions. The proof of the theorems, however, would require the following lemmas:

**LEMMA 1.** Let $\{a_T\}$ and $\{b_T\}$, $T \in [0,\infty)$, be two sequences. If $\lim_{T\to\infty} b_T = 0$, then

$$\liminf_{T\to\infty}(a_T \pm b_T) = \liminf_{T\to\infty} a_T \pm \lim_{T\to\infty} b_T = \liminf_{T\to\infty} a_T$$

$$\left(\limsup_{T\to\infty}(a_T \pm b_T) = \limsup_{T\to\infty} a_T \pm \lim_{T\to\infty} b_T = \limsup_{T\to\infty} a_T\right).$$

**LEMMA 2.** Let $\{a_T\}$ and $\{b_T\}$, $T \in [0,\infty)$, be two sequences. If $\liminf_{T\to\infty} a_T > 0$, $\lim_{T\to\infty} b_T > 0$, then $\liminf_{T\to\infty}(ab)_T = \liminf_{T\to\infty} a_T \cdot \lim_{T\to\infty} b_T$.

Proof. See the proof of Lemma 1 in Ref. [3].

**THEOREM 1.** Assume that $x^\circ(t) \equiv \lim_{T\to\infty} x_T(t)$ exist. Assume also that

$\lim_{T\to\infty} \sum_{t=0}^{T} \left(v\left(x^\circ(t), x^\circ(t+1), t\right)\right)$ is $\infty$ or a finite value. If

$\lim_{T\to\infty} \left\{ \sum_{t=0}^{T} \left\{ \left(v\left(x_T(t), x_T(t+1), t\right) - v\left(\tilde{x}^{\circ T}(t), \tilde{x}^{\circ T}(t+1), t\right)\right) \right\} \right\} = 0$, then no other infinite horizon attainable path $\{x_1(t)\}$ overtakes $\{x^\circ(t)\}$, in the sense of the modified overtaking criterion as defined in Definition 3.

Proof: For any attainable path $\{x_1(t)\}$, from Lemma 1, we see that

$$\liminf_{T\to\infty}\left\{\sum_{t=0}^{T}\left\{\left(v(x_1(t),x_1(t+1),t)-v(x_T^*(t),x_T^*(t+1),t)\right)\right\}\right\}$$

$$=\liminf_{T\to\infty}\left\{\sum_{t=0}^{T}\left\{\left(v(x_1(t),x_1(t+1),t)-v(\tilde{x}^{\circ T}(t),\tilde{x}^{\circ T}(t+1),t)\right)-\left(\left(v(x_T(t),x_T(t+1),t)-v(\tilde{x}^{\circ T}(t+1),\tilde{x}^{\circ T}(t+1),t)\right)\right)\right\}\right\}$$

$$=\liminf_{T\to\infty}\left\{\sum_{t=0}^{T}\left(v(x_1(t),x_1(t+1),t)-v(\tilde{x}^{\circ T}(t),\tilde{x}^{\circ T}(t),t)\right)\right\}-\lim_{T\to\infty}\left\{\sum_{t=0}^{T}\left(v(x_T(t),x_T(t),t)-v(\tilde{x}^{\circ T}(t),\tilde{x}^{\circ T}(t),t)\right)\right\}.$$

Since $\lim_{T\to\infty}\left\{\sum_{t=0}^{T}\left\{\left(v(x_T(t),x_T(t+1),t)-v(\tilde{x}^{\circ T}(t),\tilde{x}^{\circ T}(t),t)\right)\right\}\right\}=0$, then

$$\liminf_{T\to\infty}\left\{\sum_{t=0}^{T}\left(v(x_1(t),x_1(t+1),t)-U_t(\tilde{x}^{\circ T}(t),\tilde{x}^{\circ T}(t),t)\right)\right\}-\lim_{T\to\infty}\left\{\sum_{t=0}^{T}\left(v(x_T(t),x_T(t+1),t)-v(\tilde{x}^{\circ T}(t),\tilde{x}^{\circ T}(t+1),t)\right)\right\}$$

$$=\liminf_{T\to\infty}\left\{\sum_{t=0}^{T}\left(v(x_1(t),x_1(t+1),t)-v(\tilde{x}^{\circ T}(t),\tilde{x}^{\circ T}(t+1),t)\right)\right\}.$$

Because $\sum_{t=0}^{T}\left\{\left(v(x_1(t),x_1(t))-v(x_T^*(t),x_T^*(t),t)\right)\right\}\leq 0$, we see that

$$\liminf_{T\to\infty}\left\{\sum_{t=0}^{T}\left(v(x_1(t),x_1(t+1),t)-v(\tilde{x}^{\circ T}(t),\tilde{x}^{\circ T}(t+1),t)\right)\right\}\leq 0. \qquad Q.E.D.$$

We have thus identified the conditions under which the limit of the solutions for the finite horizon problems is optimal among all attainable paths for the infinite horizon optimization problem in $R$, in the sense of a modified overtaking criterion. Obviously, if the condition for the "agreeable plans" in the sense of Ref. [6] is satisfied (for our case, it would require $\lim_{T\to\infty}\left\{\sum_{t=0}^{T}\left\{\left(v(x_T^*(t),x_T^*(t+1),t)-v(\tilde{x}^{\circ T}(t),\tilde{x}^{\circ T}(t+1),t)\right)\right\}\right\}=0$), the condition in Theorem 1 is also satisfied. Hence, we have presented a weaker condition, which only requires that the loss accompanying the limit practice to be negligible as compared to the sum of the utility sequence to be maximized.

## 3.2 Transversality Condition: Necessity

Suppose that the optimal path to (1) exists and is given by $\{x^*(t)\}$, optimal in the sense of the overtaking criterion as defined in Definition 3, next we consider the necessity of the TVC.

We perturb the optimal path $\{x^*(t)\}$ to derive an attainable path $\{x(t)\}$:

$\exists p(t) > 0,\ p(0) = 0,\ \exists \delta > 0,\quad$ such that $\{x(t) \equiv x^*(t) + \varepsilon \cdot p(t), \varepsilon \in (-\delta, \delta)\} \in X$.

We define

$$J(\varepsilon, T) = \inf_{T \leq T'} \sum_{t=0}^{T'} \left( v\left(\tilde{x}^T(t), \tilde{x}^T(t+1), t\right) - v\left(\tilde{x}^{*T}(t), \tilde{x}^{*T}(t+1), t\right) \right)$$

$$= \inf_{T \leq T'} \sum_{t=0}^{T'-1} \left( v\left(\tilde{x}(t), \tilde{x}(t+1), t\right) - v\left(\tilde{x}^*(t), \tilde{x}^*(t+1), t\right) \right) + v\left(x(T'), 0, T'\right) - v\left(x^*(T'), 0, T'\right).$$

Note that $J(0, T) = 0$. Because the path $\{x^*(t)\}$ is optimal following Definition 3, we see that

$$\lim_{T \to \infty} J(\varepsilon, T) \equiv \liminf_{T \to \infty} \left( \sum_{t=0}^{T-1} \left( v\left(\tilde{x}(t), \tilde{x}(t+1), t\right) - v\left(\tilde{x}^*(t), \tilde{x}^*(t+1), t\right) \right) + v\left(x(T), 0, T\right) - v\left(x^*(T), 0, T\right) \right) \leq 0.$$

In general, $\dfrac{d}{d\varepsilon} \lim_{T \to 0} f(\varepsilon, T) = \lim_{T \to 0} \dfrac{d}{d\varepsilon} f(\varepsilon, T)$ only if $\lim_{T \to 0} \dfrac{d}{d\varepsilon} f(\varepsilon, T)$ converges uniformly for $\varepsilon$ (Ref. [8]). We then assume

**ASSUMPTION 1.** Assume that $\dfrac{J(\varepsilon)}{\varepsilon}$ converges uniformly for $\varepsilon$ when $T \to \infty$.

**ASSUMPTION 2.** We assume for any $T > 0$,

$$\inf_{T \leq T'} \sum_{t=0}^{T'} \left( v\left(\tilde{x}^T(t), \tilde{x}^T(t+1), t\right) - v\left(\tilde{x}^{*T}(t), \tilde{x}^{*T}(t+1), t\right) \right) \text{ converges uniformly for } \varepsilon.$$

Note that in Assumption 2, $x(t) \equiv x^*(t) + \varepsilon \cdot p(t), \varepsilon \in (-\delta, \delta)$.

**THEOREM 2.** For any interior optimal path $\{x^*(t)\}$, the TVC to (1) is given by

$$\liminf_{T \to \infty} \left( v_2\left(x^*(T-1), x_2^*(T), T-1\right) + v_1\left(x^*(T), 0, T\right) \right) x^*(T) \geq 0.$$

**Proof.** Under Assumption 1 and 2, we see that

$$\left. \frac{d}{d\varepsilon} \right|_{\varepsilon > 0} \lim_{T \to \infty} J(\varepsilon, T)$$

$$= \liminf_{T \to \infty} \sum_{t=0}^{T-1} \left( v_1\left(\tilde{x}^*(t), \tilde{x}^*(t+1), t\right) p(t) \right) + \left( v_2\left(\tilde{x}^*(t), \tilde{x}^*(t+1), t\right) p(t) \right) + v_1\left(x^*(T), 0, T\right) p(T)$$

$$= \liminf_{T \to \infty} \sum_{t=0}^{T-1} \left( v_1\left(\tilde{x}^*(t), \tilde{x}^*(t+1), t\right) \right) + \left( v_2\left(\tilde{x}^*(t), \tilde{x}^*(t+1), t\right) \right) p(t)$$

$$+ \left( v_2\left(x^*(T-1), x^*(T), T-1\right) + v_1\left(x^*(T), 0, T\right) \right) p(T) \leq 0,$$

(A)

whereas

$$\left. \frac{d}{d\varepsilon} \right|_{\varepsilon < 0} \lim_{T \to \infty} J(\varepsilon, T)$$

$$= \liminf_{T \to \infty} \sum_{t=0}^{T-1} \left( v_1\left(\tilde{x}^*(t), \tilde{x}^*(t+1), t\right) p(t) \right) + \left( v_2\left(\tilde{x}^*(t), \tilde{x}^*(t+1), t\right) p(t) \right) + v_1\left(x^*(T), 0, T\right) p(T)$$

$$= \limsup_{T \to \infty} \sum_{t=0}^{T-1} \left( v_1\left(\tilde{x}^*(t), \tilde{x}^*(t+1), t\right) \right) + \left( v_2\left(\tilde{x}^*(t), \tilde{x}^*(t+1), t\right) \right) p(t)$$

$$+ \left( v_2\left(x^*(T-1), x^*(T), T-1\right) + v_1\left(x^*(T), 0, T\right) \right) p(T) \geq 0.$$

(B)

Hence, we see that the Euler equation is then given by

$$\left( v_1\left(\tilde{x}^*(t), \tilde{x}^*(t+1), t\right) + v_2\left(\tilde{x}^*(t), \tilde{x}^{*T}(t+1), t\right) \right) p(t) = 0.$$

Because $\forall t, \exists p(t) > 0$, such that $\{x(t) \equiv x^*(t) + \varepsilon \cdot x^*(t)\} \in X$, where $p(t) \neq 0$ (For example, $p(t) = \begin{cases} p(s) = 0, s \neq t; \\ p(t) = \pm 1. \end{cases}$), we can further state the Euler equation as follows:

$$\forall t, \left(v_1\left(\tilde{x}^*(t), \tilde{x}^*(t+1), t\right) + v_2\left(\tilde{x}^*(t), \tilde{x}^{*T}(t+1), t\right)\right) = 0.$$

Furthermore, from (A) and (B), we see that the TVC is given by

$$\lim_{T \to \infty} \left(v_2\left(x^*(T-1), x^*(T), T-1\right) + v_1\left(x^*(T), 0, T\right)\right) p(T) = 0.$$

For the case where $p(t) = x^*(t)$, following (A), if we assume that $\exists \delta > 0$ such that $\{x(t) \equiv x^*(t) + \varepsilon \cdot p(t), \varepsilon \in (0, \delta)\} \in X$, the TVC is then given as

$$\liminf_{T \to \infty} \left(v_2\left(x^*(T-1), x^*(T), T-1\right) + v_1\left(x^*(T), 0, T\right)\right) x^*(T) \leq 0.$$

On the other hand, if we assume that $\exists \delta > 0$ such that $\{x(t) \equiv x^*(t) + \varepsilon \cdot p(t), \varepsilon \in (-\delta, 0)\} \in X$, following (B), we see that the TVC is then

$$\limsup_{T \to \infty} \left(v_2\left(x^*(T-1), x^*(T), T-1\right) + v_1\left(x^*(T), 0, T\right)\right) x^*(T) \geq 0,$$

which corresponds to the discussion in Kamihigashi (2001).      Q.E.D.

### 3.3    Transversality Condition: Sufficiency

We proceed to consider the sufficiency of the TVC. We impose the following assumptions on $v$:

**ASSUMPTION 3.** The function $v$ is continuous.

**ASSUMPTION 4.** The function $v$ is assumed to be non-increasing in the second argument, i.e., for all $t$, $v_{12}(\cdot) \geq 0$.

**ASSUMPTION 5.** The function $v$ is concave; i.e.

$v(\theta(x,y)+(1-\theta)(x',y')) \geq \theta v(x,y)+(1-\theta)v(x',y')$, all $(x,y),(x',y') \in X$ and all $\theta \in (0,1)$.

**ASSUMPTION 6.** The function $v$ is continuously differentiable on the interior of $X$.

Assumptions 3-6 are commonly postulated in economics and will be used to establish the sufficiency of the TVCs.

**THEOREM 3.** Let the function $v$ satisfy Assumptions 4-7. Then the path $x^*(t)$ is optimal for the problem (1), given $x_0$, if it satisfy the Euler equation given by

$$v_2(x^*(t), x^*(t+1), t) + v_1(x^*(t+1), x^*(t+2), t+1) = 0,$$

and the following TVC

$$\limsup_{T \to \infty} (v_2(x^*(T-1), x^*(T), T-1) + v_1(x^*(T), 0, T))x^*(T) \geq 0.$$

*Proof.* Let the function $v$ satisfy Assumptions 3-6. Let $x_0$ be given; let $x^*(t)$ satisfy Definition 3; and let $x(t)$ be any feasible path. It is sufficient to show that the difference, call it $D$, between the objective function in (1) evaluated at $x(t)$ and $x^*(t)$ is nonnegative:

$$D \equiv \liminf_{T \to \infty} \left\{ \sum_{t=0}^{T-1} (v(x(t), x(t+1), t) - v(x^*(t), x^*(t+1), t)) \right\} + v(x(T), 0, T) - v(x^*(T), 0, T) \leq 0.$$

Since $v$ is continuous, concave, and differentiable (Assumptions 2, 5, and 6), we have

$$\sum_{t=0}^{T-1}\left(v\left(x(t),x(t+1),t\right)-v\left(x^{*}(t),x^{*}(t+1),t\right)\right)+\left(v\left(x(T),0,T\right)-v\left(x^{*}(T),0,T\right)\right)$$

$$\leq \sum_{t=0}^{T-1}\left(v_{1}\left(x^{*}(t),x^{*}(t+1),t\right)\left(x(t)-x^{*}(t)\right)-v_{2}\left(x^{*}(t),x^{*}(t+1),t\right)\left(x(t)-x^{*}(t)\right)\right)$$

$$+v_{1}\left(x^{*}(T),0,T\right)\left(x(T)-x^{*}(T)\right).$$

Note that the last term in the above equation can be further stated as

$$\sum_{t=0}^{T-1}\left(v_{1}\left(x^{*}(t),x^{*}(t+1),t\right)\left(x(t)-x^{*}(t)\right)+v_{2}\left(x^{*}(t),x^{*}(t+1),t\right)\left(x(t)-x^{*}(t)\right)\right)$$

$$+v_{1}\left(x^{*}(T),0,T\right)\left(x(T)-x^{*}(T)\right)$$

$$=v_{1}\left(x^{*}(0),x^{*}(1),0\right)\left(x(0)-x^{*}(0)\right)+v_{2}\left(x^{*}(0),x^{*}(1),0\right)\left(x(1)-x^{*}(1)\right)$$

$$+v_{1}\left(x^{*}(1),x^{*}(2),1\right)\left(x(1)-x^{*}(1)\right)++v_{2}\left(x^{*}(1),x^{*}(2),1\right)\left(x(2)-x^{*}(2)\right)$$

$$+\cdots\cdots$$

$$+v_{1}\left(x^{*}(T-1),x^{*}(T),T-1\right)\left(x(T-1)-x^{*}(T-1)\right)+v_{2}\left(x^{*}(T-1),x^{*}(T),T-1\right)\left(x(T)-x^{*}(T)\right)$$

$$+v_{1}\left(x^{*}(T),0,T\right)\left(x(T)-x^{*}(T)\right).$$

Because $x^{*}(0)=x_{0}$, the first term in the RHS of the last equation disappears, and the above equation can be further stated as

$$\sum_{t=0}^{T-1}\left(v_{1}\left(x^{*}(t),x^{*}(t+1),t\right)\left(x(t)-x^{*}(t)\right)+v_{2}\left(x^{*}(t),x^{*}(t+1),t\right)\left(x(t)-x^{*}(t)\right)\right)$$

$$+v_{1}\left(x^{*}(T),0,T\right)\left(x(T)-x^{*}(T)\right)$$

$$=v_{2}\left(x^{*}(0),x^{*}(1),0\right)\left(x(1)-x^{*}(1)\right)+v_{1}\left(x^{*}(1),x^{*}(2),1\right)\left(x(1)-x^{*}(1)\right)$$

$$+v_{2}\left(x^{*}(1),x^{*}(2),1\right)\left(x(2)-x^{*}(2)\right)+v_{1}\left(x^{*}(2),x^{*}(3),2\right)\left(x(2)-x^{*}(2)\right)$$

$$+\cdots\cdots$$

$$+v_{2}\left(x^{*}(T-1),x^{*}(T),T-1\right)\left(x(T)-x^{*}(T)\right)+v_{1}\left(x^{*}(T),0,T\right)\left(x(T)-x^{*}(T)\right).$$

Applying the Euler's equation, i.e.

$v_{2}\left(x^{*}(t),x^{*}(t+1),t\right)+v_{1}\left(x^{*}(t+1),x^{*}(t+2),t+1\right)=0$, $\forall t$, we see that

$$\sum_{t=0}^{T-1}\left(v_1\left(x^*(t),x^*(t+1),t\right)\left(x(t)-x^*(t)\right)+v_2\left(x^*(t),x^*(t+1),t\right)\left(x(t)-x^*(t)\right)\right)$$
$$+v_1\left(x^*(T),0,T\right)\left(x(T)-x^*(T)\right)$$
$$=\left(v_2\left(x^*(T-1),x^*(T),T-1\right)+v_1\left(x^*(T),0,T\right)\right)\left(x(T)-x^*(T)\right)$$
$$=\left(v_2\left(x^*(T-1),x^*(T),T-1\right)+v_1\left(x^*(T),0,T\right)\right)x(T)$$
$$-\left(v_2\left(x^*(T-1),x^*(T),T-1\right)+v_1\left(x^*(T),0,T\right)\right)x^*(T)$$
$$=\left(-v_1\left(x^*(T),x^*(T+1),T\right)+v_1\left(x^*(T),0,T\right)\right)x(T)$$
$$-\left(-v_1\left(x^*(T),x^*(T+1),T\right)+v_1\left(x^*(T),0,T\right)\right)x^*(T)$$
$$\leq -\left(v_2\left(x^*(T-1),x^*(T),T-1\right)+v_1\left(x^*(T),0,T\right)\right)x^*(T)$$
$$\leq \liminf_{T\to\infty}\left(v_2\left(x^*(T-1),x^*(T),T-1\right)+v_1\left(x^*(T),0,T\right)\right)x^*(T)$$

where the last line uses the fact that $v_{12}(\cdot)\geq 0$ (Assumption 5) and $x(t)\geq 0$, all $t$.

Because $v_{12}(\cdot)\geq 0$ (Assumption 5) and $x(t)\geq 0$, all $t$. Hence

$$\sum_{t=0}^{T-1}\left(v\left(x(t),x(t+1),t\right)-v\left(x^*(t),x^*(t+1),t\right)\right)+\left(v\left(x(T),0,T\right)-v\left(x^*(T),0,T\right)\right)$$
$$\leq -\left(v_2\left(x^*(T-1),x^*(T),T-1\right)+v_1\left(x^*(T),0,T\right)\right)x^*(T).$$

Taking limit inf both sides of the above inequality, we derive
$$D\leq -\limsup_{T\to\infty}\left(v_2\left(x^*(T-1),x^*(T),T-1\right)+v_1\left(x^*(T),0,T\right)\right)x^*(T),$$

because $\limsup_{T\to\infty}\left(v_2\left(x^*(T-1),x^*(T),T-1\right)+v_1\left(x^*(T),0,T\right)\right)x^*(T)\geq 0,$ we see that
$D\leq 0.$ Q.E.D.

## 4. Applications

### 4.1 A One-Sector Growth Model

In this section, we first consider the applicability of the results by examining a parametric example of the one-sector growth model in economics. We consider a social planner's problem, in which $\beta \in (0,1]$ is the discount factor. The planner's objective is

(3) $$\max_{c(t)} \sum_{t=0}^{T} \beta^t \ln(c(t)), \text{ where } T \in [0, \infty),$$

subject to

(4) $$c(t) + k(t+1) = f(k(t)), \text{ where } f(k(t)) = k(t)^\alpha,$$

$$0 < \alpha < 1, \ \beta \in (0,1], \ k(0) \equiv k_0 \text{ given},$$

$$0 < k(t) < 1, \ k(T+1) = 0.$$

Note that because $0 < k(t) < 1$, we have $0 < c(t) < 1$. The optimal solution $(c_T(t), k_T(t))$ for time horizon $T$ can be found with the method of Lagrange multipliers. The optimal solution $(c_T(t), k_T(t))$ is uniquely determined given $k_0 > 0$ and $k(T+1) = 0$, with $k(t+1) = \alpha\beta \frac{1-(\alpha\beta)^{T-t}}{1-(\alpha\beta)^{T-t+1}} k^\alpha(t)$, $t = 0, 1, \ldots, T$, it is given by

$$c_T(t) = \frac{(1-(\alpha\beta))(\alpha\beta)^{\frac{\alpha(1-\alpha^t)}{1-\alpha}} k_0^{\alpha^{t+1}}}{\left(1-(\alpha\beta)^{T-t+1}\right)^{1-\alpha} \left(1-(\alpha\beta)^{T-t+2}\right)^{\alpha(1-\alpha)} \cdots \left(1-(\alpha\beta)^{T}\right)^{\alpha^{t-1}(1-\alpha)} \left(1-(\alpha\beta)^{T+1}\right)^{\alpha^t}},$$

$$k_T(t) = \frac{\left(1-(\alpha\beta)^{T-t+1}\right)(\alpha\beta)^{\frac{1-\alpha^t}{1-\alpha}} k_0^{\alpha^t}}{\left(1-(\alpha\beta)^{T-t+2}\right)^{1-\alpha} \left(1-(\alpha\beta)^{T-t+3}\right)^{\alpha(1-\alpha)} \cdots \left(1-(\alpha\beta)^{T+1}\right)^{\alpha^{t-1}}}.$$

When time approaches to infinity, their limits are

$$c^\circ(t) \equiv \lim_{T \to \infty} c_T(t) = (1-\alpha\beta)^{1/(1-\alpha)} (\alpha\beta)^{\alpha/(1-\alpha)} \left( \frac{k_0}{(\alpha\beta)^{1/(1-\alpha)}} \right)^{\alpha^{t+1}},$$

$$k^\circ(t) \equiv \lim_{T \to \infty} k_T(t) = (\alpha\beta)^{\alpha^{t-1}+\cdots+1} K_0^{\alpha^t} = (\alpha\beta)^{1/(1-\alpha)} \left( \frac{k_0}{(\alpha\beta)^{1/(1-\alpha)}} \right)^{\alpha^t},$$

$$\lambda^\circ(t) \equiv \lim_{T \to \infty} \lambda_T(t) = (1-\alpha\beta)^{-1} (\alpha\beta)^{-\alpha/(1-\alpha)} \left( \frac{k_0}{(\alpha\beta)^{1/(1-\alpha)}} \right)^{-\alpha^{t+1}},$$

where $\lambda_T(t)$ is the Lagrange multiplier for time horizon $T$, whereas $\lambda^\circ(t)$ is its limit when $T \to \infty$. It stands for the shadow price of capital at time $t$.

We first consider the TVC à la Kamihigashi, which is given by

$$\liminf_{T \to \infty} \left( -v_{T,2}(x_T^*, x_{T+1}^*) x_{T+1}^* \right) \geq 0.$$

It is easy to verify that for the current example,

$$\lim_{T \to \infty} \left( -v_{T,2}(x_T^*, x_{T+1}^*) x_{T+1}^* \right) = \frac{\alpha\beta}{1-\alpha\beta} > 0, \text{ since } \alpha \in (0,1) \text{ and } \beta \in (0,1].$$

In other words, the TVC à la Kamihigashi is not satisfied for the current model.

On the other hand, the TVC as given in Theorem 3 is

$$\liminf_{T \to \infty} \left( v_{T-1,2}(x_{T-1}^*, x_T^*) + v_{T,1}(x_T^*, 0) \right) x_T^* \leq 0,$$

showing that for the current model, we have

$$\liminf_{T \to \infty} \left( v_{T-1,2}\left(x_{T-1}^*, x_T^*\right) + v_{T,1}\left(x_T^*, 0\right) \right) x_T^*$$

$$= \liminf_{T \to \infty} \left( -\beta^{T-1} \frac{1}{\left(x_{T-1}^*\right)^\alpha - x_T^*} + \beta^T \frac{\alpha\left(x_T^*\right)^{\alpha-1}}{\left(x_T^*\right)^\alpha - 0} \right) x_T^*$$

$$= \liminf_{T \to \infty} \left( -\beta^{T-1} \frac{x_T^*}{\left(x_{T-1}^*\right)^\alpha - x_T^*} + \alpha\beta^T \right)$$

$$= \liminf_{T \to \infty} \left( \beta^{T-1} \frac{(\alpha\beta)^{1/(1-\alpha)}}{(\alpha\beta)^{1/(1-\alpha)} - \left((\alpha\beta)^{1/(1-\alpha)}\right)^\alpha} + \alpha\beta^T \right)$$

Clearly, when $\beta \in (0,1)$, $\liminf_{T \to \infty} \left( v_{T-1,2}\left(x_{T-1}^*, x_T^*\right) + v_{T,1}\left(x_T^*, 0\right) \right) x_T^* = 0$, whereas when $\beta = 1$, we see that $\liminf_{T \to \infty} \left( v_{T-1,2}\left(x_{T-1}^*, x_T^*\right) + v_{T,1}\left(x_T^*, 0\right) \right) x_T^* = \frac{1}{1-\alpha^{-1}} + \alpha < 0$ since $\alpha \in (0,1)$.

### 4.2 A Counterexample

Next we consider a simple example, the transversality condition of which takes the form of Theorem (3). However, the transversality condition á la Kamihigashi is violated. We consider the following example:

$$v \equiv -(x_t - a)^2 - b(x_{t+1} - x_t),$$

where $a > 0$, $b > 0$, and the initial value $x_0$ is given, with $x_0 = a$.

We see that Euler's equation is

$$-b - 2\left(x_{t+1}^* - a\right) + b = 0,$$

indicating $x_{t+1} = a$. Because $-v_{t,2}^* x_{T+1}^* = ba > 0$, we see that TVC á la Kamihigashi is violated. On the other hand, because

$$\left(v_{t,2} + v_{t,1}(x_t, 0)\right) x_t = (-b + b) a = 0,$$

TVC as given in Theorem 3 holds. Finally, it is also easy to show that once Assumption 1 and 2 are satisfied through proper choice of perturbation, both Theorem 1 and 2 can be satisfied.

## 5. Conclusions

In this paper, we have considered how to construct the optimal solutions for a general discrete time infinite horizon optimal control problem. We have established necessary and sufficient conditions for optimality in the sense of a modified optimality criterion. We have also considered how transversality conditions are related to steady states. We then applied our new results to two examples to demonstrate how the new transversality conditions derived in this paper differ from those given in the previous literature. Clearly, our results can be readily applied to the continuous time cases.